\documentclass{elsarticle}
\pdfoutput=1

\usepackage{mathtools, amssymb}
\usepackage{hyperref}

\mathtoolsset{showonlyrefs}

\numberwithin{equation}{section}
\newdefinition{dfn}{Definition}[section]
\newdefinition{rmk}[dfn]{Remark}
\newtheorem{lem}[dfn]{Lemma}
\newtheorem{thm}[dfn]{Theorem}
\newtheorem{cor}[dfn]{Corollary}
\newproof{prf}{Proof}
\newproof{prfMain}{Proof of Theorem~\ref{thm:main}}
\newproof{prfThm}{Proof of Theorem~\ref{thm:Landrock}}
\newproof{prfCor}{Proof of Corollary~\ref{cor:Landrock}}
\newproof{prfLem1}{Proof of Lemma~\ref{lem:adjoint}}
\newproof{prfLem2}{Proof of Lemma~\ref{lem:relation}}

\DeclareMathOperator{\rad}{rad}
\DeclareMathOperator{\soc}{soc}
\DeclareMathOperator{\Hom}{Hom}
\newcommand{\HomA}{\Hom_A}
\newcommand{\HomAop}{\Hom_{\op{A}}}
\renewcommand{\cap}{\operatorname{cap}}
\renewcommand{\mod}{\operatorname{mod}}
\newcommand{\op}[1]{{{#1}^{\mathrm{op}}}}
\newcommand{\blank}{{-}}
\newcommand{\Fdual}[1]{#1^*}
\newcommand{\Fddual}[1]{#1^{**}}
\newcommand{\Adual}[1]{#1^\vee}
\newcommand{\nakayama}[1]{\nu#1}

\begin{document}

\begin{frontmatter}

\title{A generalization of \\ dual symmetry and reciprocity \\ for symmetric algebras}
\author{Taro Sakurai}
\ead{tsakurai@math.s.chiba-u.ac.jp}
\address{Department of Mathematics and Informatics, Graduate School of Science,\\ Chiba University,\\ 1-33, Yayoi-cho, Inage-ku, Chiba-shi, Chiba, 263-8522 Japan}
\journal{Journal of Algebra}

\begin{abstract}
Slicing a module into semisimple ones is useful to study modules.
Loewy structures provide a means of doing so.
To establish the Loewy structures of projective modules over a finite dimensional symmetric algebra over a field \( F \),
the Landrock lemma is a primary tool.
The lemma and its corollary relate radical layers of projective indecomposable modules to radical layers of the \( F \)-duals of those modules (``dual symmetry'') and to socle layers of those modules (``reciprocity'').

We generalize these results to an arbitrary finite dimensional algebra \( A \).
Our main theorem,
which is the same as the Landrock lemma for finite dimensional symmetric algebras,
relates radical layers of projective indecomposable modules \( P \) to radical layers of the \( A \)-duals of those modules and to socle layers of injective indecomposable modules \( \nakayama{P} \) where \( \nakayama{(\blank)} \) is the Nakayama functor.
A key tool to prove the main theorem is a pair of adjoint functors,
which we call socle functors and capital functors.
\end{abstract}


\end{frontmatter}

\section{Introduction}
Semisimple modules are one of the most well-understood classes of modules.
Hence slicing a module into semisimple ones is a natural way to study modules.
Loewy structures provide a means of doing so.
To establish Loewy structures several studies has been done \cite{Benson, Koshitani, Waki}.
A primary tool in these studies is the Landrock lemma \cite{Landrock:J, Landrock:B},
which is stated below.

Let \( A \) be a finite dimensional algebra over a field \( F \) and
\( \Fdual{(\blank)} \coloneqq \Hom_F(\blank, F) \) the \emph{\( F \)-dual functor}.
The opposite algebra is denoted by \( \op{A} \).
The term module refers to a finitely generated right module.
Recall that \( A \) is a \emph{symmetric algebra} if \( A \cong \Fdual{A} \) as \( (A, A) \)-bimodules.
For other notations see Definition~\ref{dfn:soc-rad}.

\begin{thm}[{Landrock \cite[Theorem~B]{Landrock:J}}]
\label{thm:Landrock}
For a finite dimensional symmetric algebra \( A \) over a field \( F \),
let \( P_i \) and \( P_j \) be the projective covers of simple \( A \)-modules \( S_i \) and \( S_j \) respectively.
Then for an integer \( n \geq 1 \) we have an \( F \)-linear isomorphism
\begin{equation}
\HomA(\rad_n P_i, S_j) \cong \HomAop\big(\rad_n(\Fdual{P_j}), \Fdual{S_i}\big).
\end{equation}
\end{thm}

\begin{cor}
\label{cor:Landrock}
Under the same assumptions and notations of Theorem~\ref{thm:Landrock}, we have an \( F \)-linear isomorphism
\begin{equation}
\HomA(\rad_n P_i, S_j) \cong \HomA(S_i, \soc_n P_j).
\end{equation}
\end{cor}

Although these results are powerful as indicated in the beginning,
they are not applicable to algebras other than finite dimensional symmetric ones.
We generalize the above results to arbitrary finite dimensional algebras.
To state our main theorem we let \( \Adual{(\blank)} \coloneqq \HomA(\blank, A) \) be the \emph{\( A \)-dual functor} and \( \nakayama{(\blank)} \coloneqq \Fdual{\big(\Adual{(\blank)}\big)} \) the \emph{Nakayama functor}.

\begin{thm}
\label{thm:main}
For a finite dimensional algebra \( A \) over a field \( F \),
let \( P_i \) and \( P_j \) be the projective covers of simple \( A \)-modules \( S_i \) and \( S_j \) respectively.
Then for an integer \( n \geq 1 \) we have \( F \)-linear isomorphisms
\begin{align}
\HomA(\rad_n P_i, S_j) &\cong \HomAop\big(\rad_n(\Adual{P_j}), \Fdual{S_i}\big)
\label{eq:dual symmetry}
\intertext{and}
\HomA(\rad_n P_i, S_j) &\cong \HomA(S_i, \soc_n \nakayama{P_j}).
\label{eq:reciprocity}
\end{align}
\end{thm}

This paper is organized as follows.
In Section~\ref{sec:Proof of Main Theorem} we introduce key tools to prove our main theorem, socle functors and capital functors,
and prove Theorem~\ref{thm:main} assuming lemmas.
We also derive Theorem~\ref{thm:Landrock} and Corollary~\ref{cor:Landrock} from our main theorem.
Section~\ref{sec:Example} deals with a simple example to see how our main theorem looks in a concrete situation.
Section~\ref{sec:Proofs of Lemmas} consists of proofs of the lemmas in Section~\ref{sec:Proof of Main Theorem}.

\section{Proof of Main Theorem}
\label{sec:Proof of Main Theorem}
We introduce basic terminology of this paper and state some useful lemmas first.
Then a proof of the main theorem is given.
The Landrock lemma is proved as a special case of the main theorem.
\begin{dfn}
\label{dfn:soc-rad}
For a module \( V \) over an algebra, \( \soc V \) denotes the sum of minimal submodules of \( V \) and \( \rad V \) denotes the intersection of maximal submodules of \( V \).
For an integer \( n \geq 0 \), the \emph{\( n \)th socle} of \( V \) is defined inductively by \( \soc^0 V = 0 \) and
\begin{equation}
\soc^n V = \{\, v \in V \mid v + \soc^{n-1} V \in \soc(V/\soc^{n-1} V) \,\}
\end{equation}
if \( n > 0 \).
For an integer \( n \geq 0 \), the \emph{\( n \)th radical} of \( V \) is also defined inductively by \( \rad^0 V = V \) and
\begin{equation}
\rad^n V = \rad(\rad^{n-1} V)
\end{equation}
if \( n > 0 \).
We then write
\begin{equation}
\soc_n V = \soc^n V / \soc^{n-1} V
\end{equation}
and call it the \emph{\( n \)th socle layer} of \( V \) for \( n \geq 1 \).
We also write
\begin{equation}
\rad_n V = \rad^{n-1} V / \rad^n V
\end{equation}
and call it the \emph{\( n \)th radical layer} of \( V \) for \( n \geq 1 \).
\end{dfn}

\begin{dfn}
For an integer \( n \geq 0 \) and a module \( V \) over an algebra we write \( \cap^n V = V / \rad^n V \) and call it the \emph{\( n \)th capital} of \( V \).
Since any homomorphism maps the \( n \)th socle into the \( n \)th socle and the \( n \)th radical into the \( n \)th radical, \( \soc^n \) and \( \cap^n \) define endofunctors.
We call these endofunctors the \emph{\( n \)th socle functor} and the \emph{\( n \)th capital functor} respectively.
\end{dfn}

The next simple lemma, which does not appear in the literature to the best of author's knowledge, is vital to prove the main theorem.
The category of finitely generated right \( A \)-modules is denoted by \( \mod A \).
\begin{lem}
\label{lem:adjoint}
Let \( A \) be a finite dimensional algebra over a field.
Then for any integer \( n \geq 0 \) the \( n \)th capital functor and \( n \)th socle functor yield an adjoint pair of functors.
\begin{equation}
\raisebox{-0.4\height}{\includegraphics{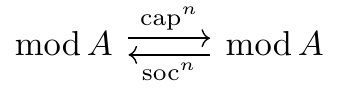}} 
\qquad
\cap^n \dashv \soc^n.
\end{equation}
\end{lem}

For a sense of unity we adopt an alias
\begin{equation}
\cap_n = \rad_n.
\end{equation}
Note that \( \cap_n \) and \( \soc_n \) define endofunctors as \( \cap^n \) and \( \soc^n \).
The following lemma can essentially be found in \cite[Problem~2.14.ii]{Nagao-Tsushima}.

\begin{lem}
\label{lem:relation}
Let \( A \) be a finite dimensional algebra over a field.
For any integer \( n \geq 1 \) we have a natural isomorphism
\begin{equation}
\raisebox{-0.5\height}{\includegraphics{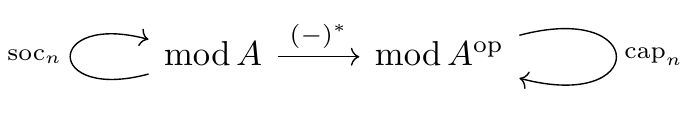}} 
\qquad
\Fdual{(\soc_n(\blank))} \cong \cap_n(\Fdual{(\blank)}).
\end{equation}
\end{lem}

\begin{prfMain}
Let us prove the reciprocity part \eqref{eq:reciprocity} first.
From the definitions we have two short exact sequences
\begin{equation}
\includegraphics{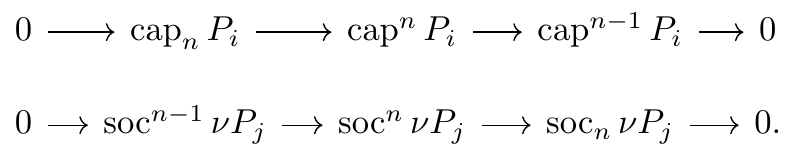}
\end{equation}
By applying exact functors \( \HomA(\blank, \nakayama{P_j}) \) and \( \HomA(P_i, \blank) \)  we have
\begin{equation}
\includegraphics{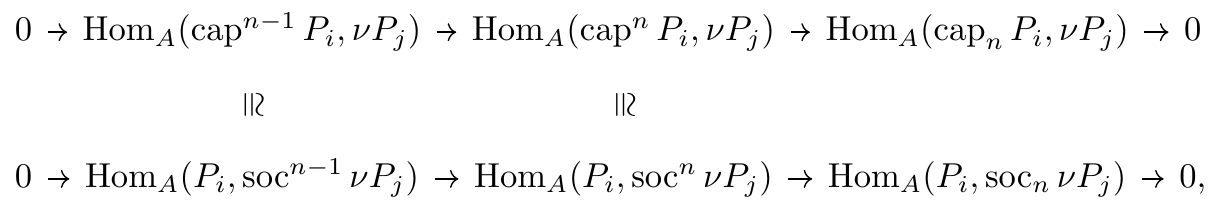}
\end{equation}
where vertical isomorphisms follow from Lemma~\ref{lem:adjoint}.
Hence we get
\begin{equation}
\HomA(\cap_n P_i, \nakayama{P_j})
\cong
\HomA(P_i, \soc_n \nakayama{P_j}). \label{eq:mid}
\end{equation}
Since the left hand side of \eqref{eq:mid} can be transformed as
\begin{align}
\HomA(\cap_n P_i, \nakayama{P_j})
&\cong
	\HomA(\cap(\cap_n P_i), \nakayama{P_j})
	&
	&
\\ &\cong
	\HomA(\cap_n P_i, \soc \nakayama{P_j})
	&
	&\text{(By Lemma~\ref{lem:adjoint}.)}
\\ &\cong
	\HomA(\cap_n P_i, S_j)
	&
	&\text{(By \cite[Lemma~III.5.1.ii]{Skowronski-Yamagata}.)}
\intertext{and the right hand side of \eqref{eq:mid} can be transformed as}
\HomA(P_i, \soc_n \nakayama{P_j})
&\cong
	\HomA(P_i, \soc(\soc_n \nakayama{P_j}))
	&
	&
\\ &\cong
	\HomA(\cap P_i, \soc_n \nakayama{P_j})
	&
	&\text{(By Lemma~\ref{lem:adjoint}.)}
\\ &\cong
	\HomA(S_i, \soc_n \nakayama{P_j}),
	&
	&
\end{align}
we have the desired isomorphism
\begin{equation}
\HomA(\cap_n P_i, S_j)
\cong
\HomA(S_i, \soc_n \nakayama{P_j}).
\end{equation}

Now let us prove the dual symmetry part \eqref{eq:dual symmetry}.
It follows immediately from the reciprocity part.
\begin{align}
\HomA(\cap_n P_i, S_j)
&\cong
\HomA(S_i, \soc_n \nakayama{P_j}) & &\text{(By \eqref{eq:reciprocity}.)}
\\ &\cong
\HomAop(\Fdual{(\soc_n \nakayama{P_j})}, \Fdual{S_i}) & &\text{(By \cite[Lemma~2.8.6.ii]{Nagao-Tsushima}.)}
\\ &\cong
\HomAop(\cap_n(\Fdual{(\nakayama{P_j})}), \Fdual{S_i}) & &\text{(By Lemma~\ref{lem:relation}.)}
\\ &\cong
\HomAop(\cap_n(\Adual{P_j}), \Fdual{S_i}).
\end{align}
\qed
\end{prfMain}

\begin{rmk}
The above proof is element-free.
\end{rmk}

We derive Theorem~\ref{thm:Landrock} and Corollary~\ref{cor:Landrock} from the main theorem in the following.
The next characterization of symmetric algebras is well-known \cite[Theorem~3.1]{Rickard}.

\begin{thm}
\label{thm:symmetric}
A finite dimensional algebra over a field is symmetric if and only if \( \Fdual{(\blank)} \cong \Adual{(\blank)} \).
\end{thm}

\begin{cor}
\label{cor:v=id}
If a finite dimensional algebra over a field is symmetric then the Nakayama functor is naturally isomorphic to the identity functor.
\end{cor}

\begin{prfThm}
Apply Theorem~\ref{thm:symmetric} to \eqref{eq:dual symmetry}.
\qed
\end{prfThm}

\begin{prfCor}
Apply Corollary~\ref{cor:v=id} to \eqref{eq:reciprocity}.
\qed
\end{prfCor}

\begin{rmk}
Okuyama and Tsushima gave a short proof of the Landrock lemma for group algebras in \cite[Theorem~2]{Okuyama-Tsushima}.
\end{rmk}

\section{Example}
\label{sec:Example}
In this section, Theorem~\ref{thm:main} is illustrated by a simple example, a connected selfinjective Naka\-yama algebra.
Let us describe this algebra and their modules first.
See \cite{Skowronski-Yamagata} for terminology.

\subsection{Algebra and Modules}
Let \( k \) and \( \ell \) be natural numbers and \( F \) an algebraically closed field. The quiver defined by
\begin{equation}
\includegraphics{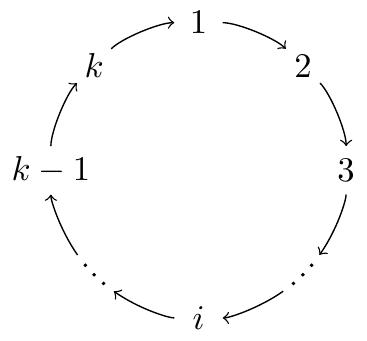}
\end{equation}
is denoted by \( \Delta_k \).
The ideal of path algebra \( F\Delta_k \) over the field \( F \) generated by the arrows of \( \Delta_k \) is denoted by \( R_k \).
Set \( N_k^\ell = F\Delta_k/{R_k}^{\ell+1} \).
This is a connected selfinjective Nakayama algebra \cite[Theorem~IV.6.15]{Skowronski-Yamagata}.
The simple module corresponding to the vertex \( i \) is denoted by \( S_i \) for \( i \in \mathbb{Z}/k\mathbb{Z} \).
Since \( N_k^\ell \) is a Nakayama algebra, the projective cover \( P_i \) of \( S_i \) and the injective envelope \( I_i \) of \( S_i \) are uniserial.
The composition series, radical series, and socle series of those modules coincide.
Those structures are described by its composition factors as below.
\begin{equation}
P_i
\quad
=
\quad
\raisebox{-0.4\height}{\includegraphics{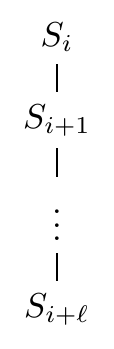}} 
\qquad
\qquad
\qquad
\qquad
I_i
\quad
=
\quad
\raisebox{-0.4\height}{\includegraphics{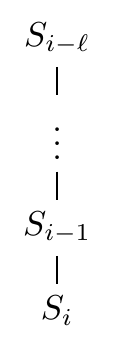}} 
\end{equation}
The Nakayama functor thus acts as \( \nakayama{P_j} = I_j = P_{j-\ell} \) for \( j \in \mathbb{Z}/k\mathbb{Z} \).

\subsection{Reciprocity}
The reciprocity \eqref{eq:reciprocity} is then illustrated by the following.
\begin{equation}
P_i
\quad
=
\quad
\raisebox{-0.5\height}{\includegraphics{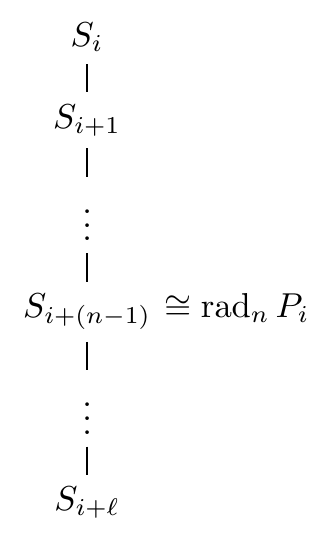}} 
\qquad
\raisebox{-0.5\height}{\includegraphics{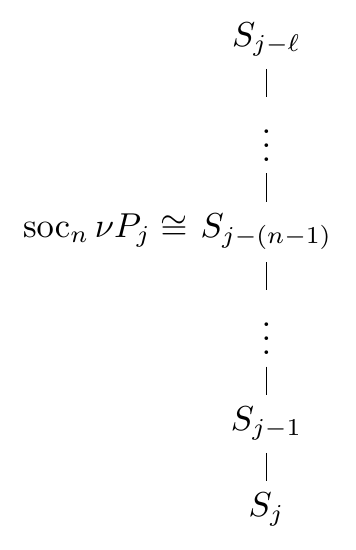}} 
\quad
=
\quad
\nakayama{P_j}
\end{equation}

\begin{equation}
\dim\HomA(\rad_n P_i, S_j) = \delta_{i+(n-1), j} = \delta_{i, j-(n-1)} = \dim\HomA(S_i, \soc_n \nakayama{P_j})
\end{equation}

\subsection{Dual Symmetry}
The radical layers of \( \Adual{P_j} \) is obtained as follows.
For \( 1 \leq n \leq \ell + 1 \) we have
\begin{align}
\rad_n(\Adual{P_j})
&\cong
\cap_n\big(\Fddual{(\Adual{P_j})}\big)
\\ &\cong
\Fdual{(\soc_n \nakayama{P_j})} & &\text{(By Lemma~\ref{lem:relation}.)}
\\ &\cong
\Fdual{S_{j-(n-1)}}.
\end{align}
The dual symmetry \eqref{eq:dual symmetry} is hence illustrated by the following.
\begin{equation}
P_i
\quad
=
\quad
\raisebox{-0.5\height}{\includegraphics{tikz/loewy3}} 
\qquad
\qquad
\raisebox{-0.5\height}{\includegraphics{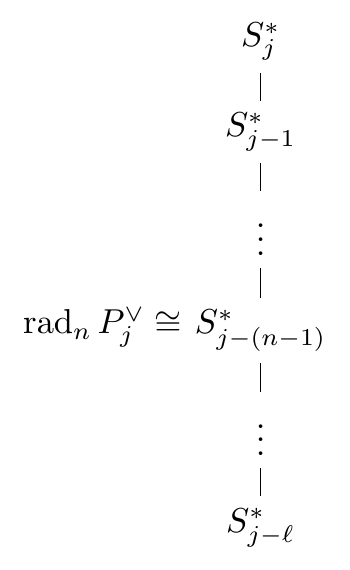}} 
\quad
=
\quad
\Adual{P_j}
\end{equation}

\begin{equation}
\dim\HomA(\rad_n P_i, S_j) = \delta_{i+(n-1), j} = \delta_{i, j-(n-1)} = \dim\HomAop(\rad_n \Adual{P_j}, \Fdual{S_i})
\end{equation}

\begin{rmk}
The examples \( N_k^\ell \) contain the symmetric cases to which the Landrock lemma and its corollary are applicable.
Indeed, the algebra \( N_k^\ell \) is symmetric if and only if \( k \) divides \( \ell \) \cite[Corollary~IV.6.16]{Skowronski-Yamagata}.
\end{rmk}

\section{Proofs of Lemmas}
\label{sec:Proofs of Lemmas}

\begin{prfLem1}
Let \( U \) and \( V \) be \( A \)-modules.
Define \( F \)-linear maps \( \eta_{U, V} \) and \( \xi_{U, V} \) by the following.
\begin{equation}
\includegraphics{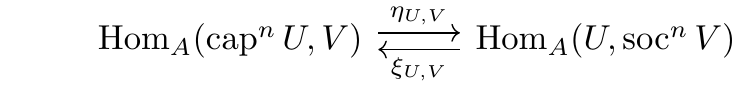}
\end{equation}
\begin{align}
&\eta_{U, V}(f) \colon u \mapsto f(u + \rad^n U)
&
&\big(f \in \HomA(\cap^n U, V)\big) \\
&\xi_{U, V}(g) \colon u + \rad^n U \mapsto g(u)
&
&\big(g \in \HomA(U, \soc^n V)\big)
\end{align}
The well-definedness of these maps follow from \( \rad^n U = U(\rad^n A) \) and \( \soc^n V = \{\, v \in V \mid v(\rad^n A) = 0 \,\} \).
It is routine work to check that these yield mutually inverse natural transformations.
\qed
\end{prfLem1}

We prove Lemma~\ref{lem:relation} in the following.
Let us introduce a lemma to prove the lemma.

\begin{lem}
Let \( A \) be a finite dimensional algebra over a field.
For any integer \( n \geq 0 \) we have a natural isomorphism
\begin{equation}
\raisebox{-0.4\height}{\includegraphics{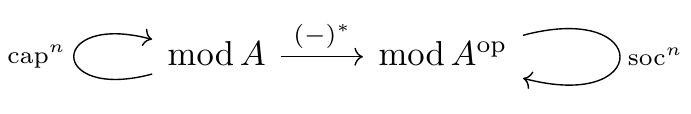}} 
\qquad
\soc^n(\Fdual{(\blank)}) \cong \Fdual{(\cap^n(\blank))}.
\end{equation}
\end{lem}

\begin{prf}
Let \( U \) be an \( A \)-module.
Then \( F \)-linear maps \( \eta_{n, U} \) and \( \xi_{n, U} \) defined by the following yield well-defined natural transformations.
\begin{equation}
\includegraphics{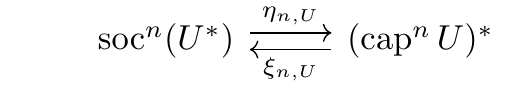}
\end{equation}
\begin{align}
&\eta_{n, U}(\lambda) \colon u +\rad^n U \mapsto \lambda(u)
&
&\big(\lambda \in \soc^n(\Fdual{U})\big) \\
&\xi_{n, U}(\mu) \colon u \mapsto \mu(u + \rad^n U)
&
&\big(\mu \in \Fdual{(\cap^n U)}\big)
\label{eq: nat trans}
\end{align}
\qed
\end{prf}

\begin{prfLem2}
To obtain a desired natural transformation, it is suffice to prove a similar statement
\begin{equation}
\soc_n(\Fdual{(\blank)}) \cong \Fdual{(\cap_n(\blank))}
\label{eq: similar}
\end{equation}
since
\begin{equation}
\Fdual{(\soc_n(\blank))}
\cong
\Fdual{(\soc_n(\Fddual{(\blank)})}
\stackrel{\text{\eqref{eq: similar}}}{\cong}
\Fddual{(\cap_n(\Fdual{(\blank)}))}
\cong
\cap_n(\Fdual{(\blank)}).
\end{equation}
For an \( A \)-module \( U \), consider the following diagram with exact rows
\begin{equation}
\includegraphics{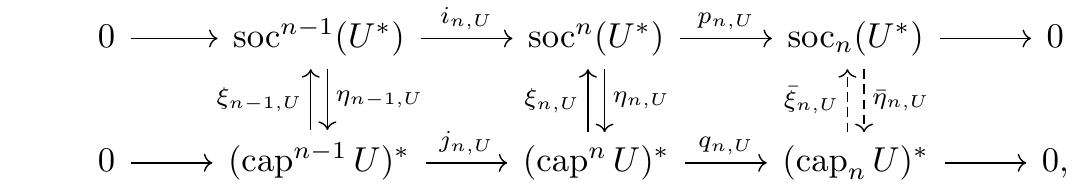}
\end{equation}
where
\begin{itemize}
\item vertical arrows are the ones defined in \eqref{eq: nat trans},
\item \( i_{n, U} \) and \( p_{n, U} \) are the inclusion and the canonical projection, and
\item \( j_{n, U} \) and \( q_{n, U} \) are homomorphisms induced from the canonical projection and the inclusion by the \( F \)-dual functor.
\end{itemize}
Since the left rectangle commutes and each rows are exact, we have unique well-defined homomorphisms \( \bar{\eta}_{n, U} \) and \( \bar{\xi}_{n, U} \) that commute the right rectangle.
A chase revels that these are mutually inverse since \( p_{n, U} \) and \( q_{n, U} \) are epimorphism.

Let \( U \) and \( V \) be \( A \)-modules and \( f \in \HomA(U, V) \).
A chase of the following diagram, which each face except the front commutes, revels that the isomorphism \( \bar{\eta}_{n, U} \) is natural since \( p_{n, V} \) is epimorphism.
\begin{equation}
\includegraphics{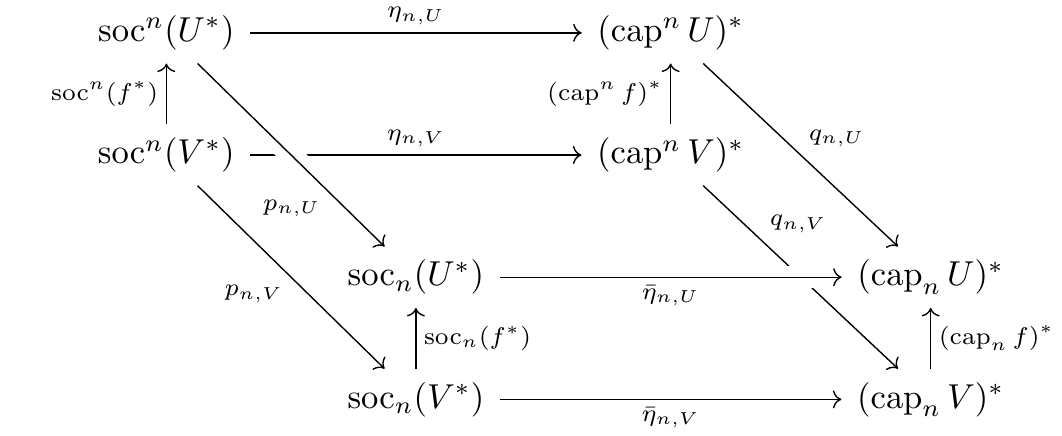}
\end{equation}
\qed
\end{prfLem2}

\section*{Acknowledgements}
The author is greatly indebted to his supervisor, Professor Shigeo Koshitani, for drawing the author's attention to the Landrock lemma and for his constant encouragement.
The author also wishes to express his thanks to Professor Osamu Iyama for providing an opportunity to deliver a talk about it and his suggestion to think more structurally.
The author would like to thank Mr.~Yoshihiro Otokita for pointing out an error in a draft.

\bibliography{sakurai}

\begin{thebibliography}{9}
\expandafter\ifx\csname natexlab\endcsname\relax\def\natexlab#1{#1}\fi
\expandafter\ifx\csname url\endcsname\relax
  \def\url#1{\texttt{#1}}\fi
\expandafter\ifx\csname urlprefix\endcsname\relax\def\urlprefix{URL }\fi

\bibitem[{Benson(1983)}]{Benson}
Benson, D., 1983. The {L}oewy structure of the projective indecomposable
  modules for {$A_{8}$} in characteristic {$2$}. Comm. Algebra 11~(13),
  1395--1432.
\newline\urlprefix\url{http://dx.doi.org/10.1080/00927878308822912}

\bibitem[{Koshitani(1985)}]{Koshitani}
Koshitani, S., 1985. On the {L}oewy series of the group algebra of a finite
  {$p$}-solvable group with {$p$}-length {$>1$}. Comm. Algebra 13~(10),
  2175--2198.
\newline\urlprefix\url{http://dx.doi.org/10.1080/00927878508823271}

\bibitem[{Landrock(1983{\natexlab{a}})}]{Landrock:J}
Landrock, P., 1983{\natexlab{a}}. The {C}artan matrix of a group algebra modulo
  any power of its radical. Proc. Amer. Math. Soc. 88~(2), 205--206.
\newline\urlprefix\url{http://dx.doi.org/10.2307/2044699}

\bibitem[{Landrock(1983{\natexlab{b}})}]{Landrock:B}
Landrock, P., 1983{\natexlab{b}}. Finite group algebras and their modules.
  Vol.~84 of London Mathematical Society Lecture Note Series. Cambridge
  University Press, Cambridge.
\newline\urlprefix\url{http://dx.doi.org/10.1017/CBO9781107325524}

\bibitem[{Nagao and Tsushima(1989)}]{Nagao-Tsushima}
Nagao, H., Tsushima, Y., 1989. Representations of finite groups. Academic
  Press, Inc., Boston, MA, translated from the Japanese.

\bibitem[{Okuyama and Tsushima(1986)}]{Okuyama-Tsushima}
Okuyama, T., Tsushima, Y., 1986. On a conjecture of {P}. {L}androck. J. Algebra
  104~(1), 203--208.
\newline\urlprefix\url{http://dx.doi.org/10.1016/0021-8693(86)90247-4}

\bibitem[{Rickard(2002)}]{Rickard}
Rickard, J., 2002. Equivalences of derived categories for symmetric algebras.
  J. Algebra 257~(2), 460--481.
\newline\urlprefix\url{http://dx.doi.org/10.1016/S0021-8693(02)00520-3}

\bibitem[{Skowro{\'n}ski and Yamagata(2011)}]{Skowronski-Yamagata}
Skowro{\'n}ski, A., Yamagata, K., 2011. Frobenius algebras. {I}. EMS Textbooks
  in Mathematics. European Mathematical Society (EMS), Z\"urich, basic
  representation theory.
\newline\urlprefix\url{http://dx.doi.org/10.4171/102}

\bibitem[{Waki(1993)}]{Waki}
Waki, K., 1993. The {L}oewy structure of the projective indecomposable modules
  for the {M}athieu groups in characteristic {$3$}. Comm. Algebra 21~(5),
  1457--1485.
\newline\urlprefix\url{http://dx.doi.org/10.1080/00927879308824631}

\end{thebibliography}
\end{document}